\documentclass[11pt]{amsart}

\pdfoutput=1

\usepackage[text={400pt,660pt},centering]{geometry}
\usepackage[foot]{amsaddr}

\usepackage{amsmath,amssymb}
\usepackage{esint,amssymb} 
\usepackage{graphicx}
\usepackage{MnSymbol}
\usepackage{mathtools} 
\usepackage[colorlinks=true, pdfstartview=FitV, linkcolor=blue, citecolor=blue, urlcolor=blue,pagebackref=false]{hyperref}
\usepackage{microtype}

\usepackage{bm}
\usepackage{dsfont}

\newtheorem{proposition}{Proposition}

\theoremstyle{remark}
\newtheorem{remark}[proposition]{Remark}

\theoremstyle{definition}
\newtheorem{definition}[proposition]{Definition}

\numberwithin{equation}{section}
\numberwithin{proposition}{section}
\numberwithin{figure}{section}
\numberwithin{table}{section}

\newcommand{\R}{\mathbb{R}}

\newcommand{\E}{\mathbb{E}}

\renewcommand{\le}{\leqslant}
\renewcommand{\ge}{\geqslant}

\makeatletter
\newsavebox\myboxA
\newsavebox\myboxB
\newlength\mylenA
\newcommand*\mybar[2][0.75]{%
    \sbox{\myboxA}{$\m@th#2$}%
    \setbox\myboxB\null
    \ht\myboxB=\ht\myboxA%
    \dp\myboxB=\dp\myboxA%
    \wd\myboxB=#1\wd\myboxA
    \sbox\myboxB{$\m@th\overline{\copy\myboxB}$}
    \setlength\mylenA{\the\wd\myboxA}
    \addtolength\mylenA{-\the\wd\myboxB}%
    \ifdim\wd\myboxB<\wd\myboxA%
       \rlap{\hskip 0.5\mylenA\usebox\myboxB}{\usebox\myboxA}%
    \else
        \hskip -0.5\mylenA\rlap{\usebox\myboxA}{\hskip 0.5\mylenA\usebox\myboxB}%
    \fi}
\makeatother

\def\build#1_#2^#3{\mathrel{\mathop{\kern 0pt#1}\limits_{#2}^{#3}}}

\begin{document}

\author[L. Chevillard]{Laurent Chevillard}
\address[L. Chevillard, M. Lagoin and S.G. Roux]{Univ Lyon, Ens de Lyon, Univ Claude Bernard, CNRS, Laboratoire de Physique, 46 all\'ee d’Italie F-69342 Lyon, France}
\email{laurent.chevillard@ens-lyon.fr}

\author[M. Lagoin]{Marc Lagoin}

\author[S.G. Roux]{St\'ephane G. Roux}

\keywords{Fractional Brownian Motion, Ornstein-Uhlenbeck processes, Multiplicative Chaos}
\subjclass[2010]{60G22, 60H05}
\date{\today}

\title{Multifractal Fractional Ornstein-Uhlenbeck Processes}

\begin{abstract}
The Ornstein-Uhlenbeck process can be seen as a paradigm of a finite-variance and statistically stationary rough random walk. Furthermore, it is defined as the unique solution of a Markovian stochastic dynamics and shares the same local regularity as the one of the Brownian motion. A natural generalization of this process able to reproduce the local regularity of a fractional Brownian motion of parameter $H$ is provided by the fractional Ornstein-Uhlenbeck process. Based on previous works, we propose to include some Multifractal corrections to this picture using a Gaussian Multiplicative Chaos. The aforementioned process, called a Multifractal fractional Ornstein-Uhlenbeck process, is a statistically stationary finite-variance process. Its underlying dynamics is non-Markovian, although non-anticipating and causal. The numerical scheme and theoretical approach are based on a regularization procedure, that gives a meaning to this dynamical evolution, which unique solution converges towards a well-behaved stochastic process.

\end{abstract}

\maketitle

%
%
%
%
%
%

\section{Introduction}
\label{s.intro}

In many situations, for example some physical phenomena governed by thermal fluctuations such as the position of small particles suspended in a quiescent fluid, the standard Brownian motion, or Wiener process, $W(t)$  is called to model and give account of certain quantities of interest. As often required by the nature of the situation, the phenomenon itself is statistically stationary and of finite-variance, therefore calling for  the Ornstein-Uhlenbeck (OU) process. It is defined as the unique solution of the stochastic dynamics given by
\begin{align}\label{eq:OUIntro}
dX(t)=-\frac{1}{T}X(t)dt+dW(t),
\end{align}
where the increment $dW(t)$ of the Wiener process $W$ can be viewed as a Gaussian white noise, and which meaning is clear as a random distribution. The large time scale $T>0$ governs the correlation duration of the process. As a continuous-time stochastic process, $X(t)$ is nowhere differentiable, and we can say that sample paths of $X$, as those of $W$,  are H$\ddot{\text{o}}$lder continuous of any order strictly less than $H=1/2$, where $H$ is the Hurst exponent.

Other physical phenomena, for example fluid turbulence \cite{Kol41}, or some time series encountered in finance \cite{Man13,BarShe01}, exhibit in a statistically averaged sense a different Hurst exponent $H\in]0,1[$, a local regularity that can be reproduced by fractional Brownian motions \cite{ManVan68}. Once again, as required for instance by the physics of turbulence which suggests a roughness of order $H\approx 1/3$ \cite{Fri95}, such phenomena are eventually statistically stationary and are of finite variance. In this case, using fractional Ornstein-Uhlenbeck (fOU) processes \cite{CheKaw03} to model them seems more appealing. They could be defined as the unique solution of the following dynamics
\begin{align}\label{eq:fOUIntro}
dX_H(t)=-\frac{1}{T}X_H(t)dt+\text{``}dW_H(t)\text{''},
\end{align}
where the random noise $dW_H$ entering in  \ref{eq:fOUIntro} is indicated between quotes to remind that its meaning has to be clarified. Whereas fractional motions can be written as stochastic integrals in a natural way \cite{PipTaq00,Nua03,TudVie07} for any  $H\in]0,1[$, the meaning pathwise of the random noise $dW_H(t)$, and consequently the meaning of the dynamics of $X_H$ (\ref{eq:fOUIntro}), is not obvious when $H<1/2$. As proposed in \cite{Che17}, a regularization procedure over a small timescale, $\epsilon>0$, gives both a precise and explicit meaning to this random measure $dW_{\epsilon,H}(t)$. Accordingly, the unique solution $X_{\epsilon,H}(t)$ of  \ref{eq:fOUIntro} once forced by such a random noise is shown in  \cite{Che17} to converge in law, as $\epsilon\to 0$, towards the canonical fOU process of \cite{CheKaw03}. Actually, the same regularized version of the noise $dW_H(t)$,  as it is studied in \cite{Che17} and that we are using in this communication, was previously proposed in the theoretical setup developed in \cite{AloMaz00} (see Remark \ref{Rem:Alos}). 

As we can see, this regularization procedure over $\epsilon$ is necessary, when $H<1/2$, to define pathwise the dynamics schematically presented in   \ref{eq:fOUIntro}. This regularization turns out to be crucial when multifractal corrections, usually referred as the intermittency phenomenon in the literature of turbulence \cite{Fri95}, are added to the picture while considering a Gaussian multiplicative chaos \cite{Man72,Kah85,RhoVar14}. It can be formally obtained as
\begin{align}\label{eq:MCIntroFormal}
 M_{\gamma}(t)\; \text{``} = \text{''} \; e^{\gamma\widetilde{X}(t)},
\end{align}
where the Gaussian field $\widetilde{X}$ entering in  \ref{eq:MCIntroFormal} is asked to be logarithmically correlated, and $\gamma\in\R$ is a free parameter. Once again, the equality given in  \ref{eq:MCIntroFormal} is put in quotes because the meaning of the exponential of such a field, necessarily being of infinite variance, has yet to be clarified. In this article, we will show how to give a precise meaning to $M_{\gamma}$ ( \ref{eq:MCIntroFormal}) while still relying on a regularization procedure over $\epsilon$, as it is usually done \cite{RhoVar14}. Moreover, we will see that the Gaussian process $\widetilde{X}$ can actually be defined as a fOU process, solution of the stochastic dynamics provided in  \ref{eq:fOUIntro}, for a vanishing Hurst exponent $H=0$. This approach was already proposed and used in different setups to model some aspects of the random nature of fluid turbulence \cite{Che17,PerGar16,PerMor18,VigFri20,DurLes20}. This method of construction differs from other ones found in the literature \cite{FyoKho16,NeuRos18,HagNeu20}. In particular, it includes the notion of causality, and its respective stochastic evolution, to this logarithmically correlated Gaussian field.

Finally, we propose a way to incorpore into the dynamics of the fOU process (\ref{eq:fOUIntro}) a Multiplicative Chaos (\ref{eq:MCIntroFormal}), and define respectively the stochastic evolution
\begin{align}\label{eq:MfOUIntro}
dX_{H,\gamma}(t)=-\frac{1}{T}X_{H,\gamma}(t)dt+\text{``}dW_{H,\gamma}(t)\text{''},
\end{align}
which eventually makes sense only up to a regularization scale $\epsilon>0$. Taking the limit $\epsilon\to 0$, the asymptotic process obtained while considering the unique solution of \ref{eq:fOUIntro}, will be referred as the Multifractal fractional Ornstein-Uhlenbeck (MfOU) process. It can be seen as a generalization of the Multifractal random walk of \cite{BacDel01} to a causal framework, and for any $H\in]0,1[$ instead of the unique value $H=1/2$.

\smallskip
\subsection*{Organization of the paper.}
In Section \ref{Sec:Setup}, we present the general setup of the article, and especially the regularized version of the stochastic evolution proposed in \ref{eq:MfOUIntro}. We include several remarks on the various ingredients entering in the construction of this dynamics. Section \ref{Sec:StatResults} is devoted to the statement of our results. In particular, we list in Proposition \ref{Prop:MfOU} some important marginals of the MfOU process $X_{H,\gamma}$, including variance and high order moments of its increments. In Section \ref{Sec:Nums}, we propose a method to synthesize a statistically stationary solution of the dynamics (\ref{eq:MfOUIntro}) in a periodical framework, and compare our theoretical predictions to our numerical estimations. We gather the proofs of Proposition \ref{Prop:MfOU} in Section \ref{sec:Proofs}.

\section{Setup, notations and remarks}\label{Sec:Setup}

\begin{definition}\label{Def:MFOU}Consider a Wiener processes $W(t)$, note by $dW(t)$ its increments over $dt$ at time $t$. In the sequel, we will be studying the unique statistically stationary solution $X_{\epsilon,H,\gamma}(t)$ of the following linear stochastic equation
\begin{equation}\label{eq:SDEintro}
dX_{\epsilon,H,\gamma}(t)=-\frac{1}{T} X_{\epsilon,H,\gamma}(t)dt + dW_{\epsilon,H,\gamma}(t),
\end{equation}
where enters a large timescale $T>0$ and a small regularizing timescale $\epsilon>0$. For a given Hurst exponent $H\in]0,1[$, and an additional free parameter $\gamma\in\R$, the dynamical evolution proposed in \ref{eq:SDEintro} is maintained by a random noise $dW_{\epsilon,H,\gamma}(t)$ that reads
\begin{equation}\label{eq:DefdWRandom}
dW_{\epsilon,H,\gamma}(t) = \omega_{\epsilon,H,\gamma}(t)dt + \epsilon^{H-\frac{1}{2}}M_{\epsilon,\gamma}(t)dW(t),
\end{equation}
where we have introduced the statistically stationary random field
\begin{equation}\label{eq:OmegaHGammaIntro}
\omega_{\epsilon,H,\gamma}(t) = \left(H-\frac{1}{2}\right)\int_{-\infty}^t(t-s+\epsilon)^{H-\frac{3}{2}}M_{\epsilon,\gamma}(s)dW(s).
\end{equation}
The random noise $dW_{\epsilon,H,\gamma}$ (\ref{eq:DefdWRandom}) and field $\omega_{\epsilon,H,\gamma}$ (\ref{eq:OmegaHGammaIntro}) involve respectively a local version of and a linear operation on the white noise $dW$ weighted by a positive random field $M_{\epsilon,\gamma}(t)$ given by
\begin{equation}\label{eq:MCIntro}
M_{\epsilon,\gamma}(t) = e^{\gamma \widetilde{X}_{\epsilon}(t)-\gamma^2\E\left( \widetilde{X}_{\epsilon}^2\right)}.
\end{equation}
Finally, the zero-average statistically stationary Gaussian process $\widetilde{X}_{\epsilon}$ that enters in the definition of the positive random weight $M_{\epsilon,\gamma}$ (\ref{eq:MCIntro}) is independent of the white noise $dW(t)$ and is of finite variance for any $\epsilon>0$. We furthermore ask, in the limit $\epsilon\to 0$, its covariance to diverge in a logarithmic fashion at small arguments, such that
\begin{equation}\label{eq:FormalCovXTildeSetup}
\lim_{\epsilon\to 0}\E\left( \widetilde{X}_{\epsilon}(t)\widetilde{X}_{\epsilon}(t+\tau)\right)=\ln_+\left(\frac{T}{|\tau|} \right) + g(\tau),
\end{equation}
with $\ln_+ (s) = \max(\ln s, 0)$ and $g$ a bounded and continuous function of its argument which goes to 0 at large arguments. A natural and causal way to achieve this is to define $\widetilde{X}_{\epsilon}$ as the unique statistically stationary solution of an independent linear and Gaussian evolution which is based on the dynamics provided in \ref{eq:SDEintro} with $\gamma=0$, using furthermore the boundary value $H=0$. Defined such a way, $\widetilde{X}_{\epsilon}$ can be seen as a regularized fractional Ornstein-Uhlenbeck of vanishing Hurst exponent.
\end{definition}
\smallskip

\begin{remark} \label{Rem:Alos}For $\gamma=0$ and thus $M_{\epsilon,0}(t) =1$, the noise $dW_{\epsilon,H,0}(t)$ (\ref{eq:DefdWRandom}) is a well-posed Gaussian and causal process. Using the language developed in \cite{AloMaz00}, it is the stochastic differential, in the sense of It\^o, of a continuous semimartingale $W_{\epsilon,H,0}(t)$ which is a natural regularized form a fractional Brownian motion. \end{remark}

\smallskip

\begin{remark} \label{Rem:fOU} As we will see, the parameter $\gamma$ eventually governs entirely the level of non Gaussianity and multifractality of the process $X_{\epsilon,H,\gamma}(t)$. To this regard, for $\gamma\ne 0$, the dynamics provided in \ref{eq:SDEintro} can be seen as a generalization of the one of the regularized fractional Ornstein-Uhlenbeck process detailed in \cite{Che17}. Using current notation, the unique statistically stationary solution of \ref{eq:SDEintro} with $\gamma=0$ can be conveniently written as
\begin{equation}\label{eq:fOUProcInt}
X_{\epsilon,H,0}(t) = \int_{-\infty}^t e^{-\frac{t-s}{T}}dW_{\epsilon,H,0}(s).
\end{equation}
It is a finite variance Gaussian process for any $H\in]0,1[$ and finite $\epsilon>0$, and it remains so as $\epsilon\to 0$. We note by $X_{H,0}$ the Gaussian stochastic process which marginals coincide with the limiting values as $\epsilon\to 0$ of those of $X_{\epsilon,H,0}$, and we call it a fractional Ornstein-Uhlenbeck (fOU) process  \cite{CheKaw03}. We have
\begin{equation}\label{eq:VarfOU}
\E\left( X_{H,0}^2\right)\equiv \lim_{\epsilon\to 0} \E\left( X_{\epsilon,H,0}^2\right)= \frac{T^{2H}\left[ \Gamma\left( H+\frac{1}{2}\right)\right]^2}{2\sin(\pi H)},
\end{equation}
where enters the function $\Gamma(z) = \int_0^\infty x^{z-1}e^{-x}dx$. More generally, since $X_{\epsilon,H,0}$ is a zero average Gaussian process, it is fully characterized by its covariance function which can be explicitly derived and behaves in the limit of vanishing regularizing scale $\epsilon\to 0$ as, considering without loss of generality $\tau\ge 0$,
\begin{align}
&\E\left( X_{H,0}(t)X_{H,0}(t+\tau)\right)\equiv \lim_{\epsilon\to 0} \E\left( X_{\epsilon,H,0}(t)X_{\epsilon,H,0}(t+\tau)\right)\notag \\
&=\frac{1}{2\sin(\pi H)}\frac{\left[ \Gamma\left( H+\frac{1}{2}\right)\right]^2}{\Gamma\left( 2H\right)}\int_0^{\infty}\left[\frac{1}{2}e^{-\frac{\tau+h}{T}}-\frac{\tau-h}{2|\tau-h|}e^{-\frac{|\tau-h|}{T}}\right]h^{2H-1}dh,\label{eq:RemExprCorrFOU1}\\
&=T^2\left[ \Gamma\left( H+\frac{1}{2}\right)\right]^2\int_{\R}e^{2i\pi \omega \tau} \frac{|2\pi\omega|^{1-2H}}{1+4\pi^2\omega^2T^2}d\omega.\label{eq:RemExprCorrFOU2}
\end{align}
Notice that the expression of the correlation function (\ref{eq:RemExprCorrFOU1}), derived in \cite{Che17}, is especially well suited to study its decay towards 0 at large arguments, in an alternative way as it is proposed in \cite{CheKaw03}, whereas its spectral form (\ref{eq:RemExprCorrFOU2}), that can be derived and justified using more general expressions of \cite{Che17}, clarifies that indeed, the asymptotic process $X_{H,0}$ coincides with the fractional Ornstein-Uhlenbeck process of \cite{CheKaw03}, and differs from the so-called tempered fractional Brownian motion of \cite{BonDid19}.

Finally, as a consequence of \ref{eq:RemExprCorrFOU1}, it is straightforward to see that fractional Brownian motions \cite{ManVan68} and fractional Ornstein-Uhlenbeck processes  \cite{CheKaw03} share the same local regularity, as it is quantified by the second moment of the increments. Defining the increment as $\delta_\tau X_{\epsilon,H,0}(t)\equiv X_{\epsilon,H,0}(t+\tau)-X_{\epsilon,H,0}(t)$, we have
\begin{align}
\E\left( \delta_\tau X_{H,0}\right)^2&\equiv \lim_{\epsilon\to 0} \E\left( \delta_\tau X_{\epsilon,H,0}\right)^2 \label{eq:S2FOU}\\
&\build{\sim}_{\tau\to 0}^{} \frac{1}{\sin(\pi H)}\frac{\left[ \Gamma\left( H+\frac{1}{2}\right)\right]^2}{\Gamma\left( 2H+1\right)}|\tau|^{2H}.\notag
\end{align}
\end{remark}
\smallskip

\begin{remark} Let us now comment on the Gaussian field $\widetilde{X}_{\epsilon}$
entering in the construction of the regularized Multifractal fractional Ornstein-Uhlenbeck process defined in \ref{Def:MFOU}. As already mentioned, $\widetilde{X}_{\epsilon}$ can be viewed as a regularized fOU process, as it is recalled in remark \ref{Rem:fOU}, with the boundary value $H=0$. More precisely, it can be written as the following stochastic integral
\begin{equation}\label{eq:IntegralXtildeDef}
\widetilde{X}_{\epsilon}(t) = \int_{-\infty}^t e^{-\frac{t-s}{T}}d\widetilde{W}_{\epsilon,0,0}(s),
\end{equation}
where
\begin{equation}
d\widetilde{W}_{\epsilon,0,0}(t) = -\frac{dt}{2}\int_{-\infty}^t(t-s+\epsilon)^{-\frac{3}{2}}d\widetilde{W}(s) + \epsilon^{-\frac{1}{2}}d\widetilde{W}(t),
\end{equation}
where the white noise $d\widetilde{W}(t)$ is chosen independently of $dW(t)$. As we will see, the Gaussian process $\widetilde{X}_{\epsilon}$ eventually shares similar statistical properties as those of the processes constructed in Refs. \cite{FyoKho16,NeuRos18} in the asymptotic limit of vanishing $\epsilon$, although we recall that it is moreover defined as a stationary solution of a well-posed stochastic evolution.

In this case, the finiteness of variance is not guaranteed as the regularizing scale $\epsilon$ goes to zero, and indeed, it is shown in Refs. \cite{Che17,PerMor18,ApoMor20} that instead,
\begin{equation}\label{eq:DivVarXtilde}
 \E\left( \widetilde{X}_{\epsilon}^2\right)\equiv \E\left( X_{\epsilon,0,0}^2\right)\build{\sim}_{\epsilon\to 0}^{}\ln\left(\frac{1}{\epsilon}\right).
\end{equation}
Whereas the variance of $\widetilde{X}_{\epsilon}$ (\ref{eq:DivVarXtilde}) diverges logarithmically with the regularizing parameter $\epsilon$, its covariance remains bounded, and we can write, as a direct application of Eqs. \ref{eq:RemExprCorrFOU1} and \ref{eq:RemExprCorrFOU2}, for any $\tau>0$,
\begin{align}
\E\left( \widetilde{X}(t)\widetilde{X}(t+\tau)\right)&\equiv \lim_{\epsilon\to 0} \E\left( X_{\epsilon,0,0}(t)X_{\epsilon,0,0}(t+\tau)\right)\notag \\
&=\int_0^{\infty}\left[\frac{1}{2}e^{-\frac{\tau+h}{T}}-\frac{\tau-h}{2|\tau-h|}e^{-\frac{|\tau-h|}{T}}\right]h^{-1}dh,\label{eq:RemExprCorrFOUH01}\\
&=\pi T^2\int_{\R}e^{2i\pi \omega \tau} \frac{|2\pi\omega|}{1+4\pi^2\omega^2T^2}d\omega.\label{eq:RemExprCorrFOUH02}
\end{align}
Starting from \ref{eq:RemExprCorrFOUH01}, we recover that, once the limit $\epsilon\to 0$ has been taken, the variance of $\widetilde{X}$ is infinite while noticing that
\begin{equation}\label{eq:DivCovXtilde}
\E\left( \widetilde{X}(t)\widetilde{X}(t+\tau)\right)\build{\sim}_{\tau\to 0}^{}\ln\left(\frac{1}{|\tau|} \right).
\end{equation}
Because of the logarithmic behavior of the covariance function of $\widetilde{X}$  for small arguments (\ref{eq:DivCovXtilde}) and its decay towards 0 at large arguments $\tau\gtrapprox T$, we will find convenient, and useful for subsequent calculations, to write in a formal way
\begin{equation}\label{eq:FormalCovXTilde}
\E\left( \widetilde{X}(t)\widetilde{X}(t+\tau)\right)=\ln_+\left(\frac{T}{|\tau|} \right) + g(\tau),
\end{equation}
with $\ln_+ (s) = \max(\ln s, 0)$ and $g$ a continuous, even and bounded function of its argument, as it is done in \cite{RhoVar14}. We could get an exact expression for the bounded function $g$ entering in \ref{eq:FormalCovXTilde} using in particular \ref{eq:RemExprCorrFOUH01}, although, as we will see,  it will turn out to be sufficient to specify its value at the origin to understand with precision the statistical properties of the MfOU process at small scales. We find
\begin{equation}\label{eq:G(0)}
g(0)=\int_{0}^{\infty}\ln(h) e^{-h}dh\approx -0.577216,
\end{equation}
which is known as minus the Euler-Mascheroni constant.
\end{remark}

\smallskip

\begin{remark}\label{rem:DefMC} The last ingredient entering in the construction of the MfOU is the positive random weight $M_{\epsilon,\gamma}$ defined in \ref{eq:MCIntro}. It ressembles a regularized version of a Gaussian Multiplicative Chaos, as it is called in the probability theory literature \cite{RhoVar14}, although, the way it is renormalized makes it converging locally as $\epsilon \to 0$ to a random distribution, which is zero-average and unit-variance. In the best situations, where subsequent integrals involve a deterministic function $f$ that ensures dominated convergence, we will write, making use of the expression of the covariance of $\widetilde{X}$ (\ref{eq:FormalCovXTilde}),
\begin{align}\label{eq:ExpandChaosAgainstF}
\lim_{\epsilon\to 0}\;\E\;\bigg( \int_{\R}  f(t) M_{\epsilon,\gamma}(t) &dW(t)\bigg)^{2n} =\\
&\frac{(2n)!}{2^n n!}T^{2\gamma^2n(n-1)}\int_{\R^{n}}\prod_{i=1}^n f^2(t_i)\prod_{i<j=1}^n  \frac{e^{4\gamma^2 g(t_i-t_j)}}{|t_i-t_j|_+^{4\gamma^2}} \prod_{i=1}^n dt_i, \notag
\end{align}
where $|.|_+=\exp(\ln_+|.|)$ and for a finite range of values of $\gamma^2$, depending on the regularity of $f$ in the vicinity of the origin and on the requirement of the finiteness of this moment. As we will see, the very precise behavior of the random noise $M_{\epsilon,\gamma}(t) dW(t)$ when integrated against test functions, such as it is depicted in \ref{eq:ExpandChaosAgainstF}, is eventually responsible of Multifractal corrections to the monofractal Gaussian process fOU.

Let us mention that defining a Multiplicative Chaos (\ref{eq:MCIntro}) in such a causal framework is far from anecdotal. In particular, causality is crucial when modeling several aspects of fluid turbulence \cite{Che17,PerMor18}. Making use of It\^o's lemma, the regularization procedure allows to derive a stochastic evolution for the random weight  $M_{\epsilon,\gamma}$ (\ref{eq:MCIntro}) itself, and we obtain
\begin{align}\label{eq:SdeChaos}
\frac{dM_{\epsilon,\gamma}(t)}{M_{\epsilon,\gamma}(t)}=-\frac{1}{T}\left[\ln M_{\epsilon,\gamma}(t)+\gamma^2  \E\left( \widetilde{X}_{\epsilon}^2\right) +
\gamma \widetilde{\omega}_{\epsilon,0,0}(t)+\frac{\gamma^2}{2\epsilon}\right]dt+\gamma\epsilon^{-\frac{1}{2}}d\widetilde{W}(t),
\end{align}
where the causal field $\widetilde{\omega}_{\epsilon,0,0}(t)$ corresponds to the expression of  \ref{eq:OmegaHGammaIntro}, with $H=0$, $\gamma=0$, and using the increment of the independent Wiener process $\widetilde{W}$ instead of $W$. It is interesting to observe that each of the terms entering in the RHS of \ref{eq:SdeChaos} diverges pathwise as $\epsilon\to 0$, whereas the variance of the stationary solution remains bounded, with $\E M^2_{\epsilon,\gamma}=1$. Notice also that the dynamical evolution of $M_{\epsilon,\gamma}$ (\ref{eq:SdeChaos}) can be seen as a stochastic Volterra system \cite{JacPan20}.

Finally, let us remark that the random weight $M_{\epsilon,\gamma}$ (\ref{eq:MCIntro}), which is a log-normal random field since it is taken as the exponential of a Gaussian process, can be seen as a particular case of the more general class of log-infinitely divisible measures \cite{BarMan02,SchMar01,BacMuz03,ChaRie05,RhoSoh14,RhoVar14}. 
\end{remark}

\smallskip 

\section{Statement of the results} \label{Sec:StatResults}
\begin{proposition}\label{Prop:MfOU}
The unique solution of the stochastic differential equation \ref{eq:SDEintro} goes at large time towards a statistically stationary process $X_{\epsilon,H,\gamma}(t)$ that we can write
\begin{equation}\label{eq:PropMfOU}
X_{\epsilon,H,\gamma}(t) = \int_{-\infty}^t e^{-\frac{t-s}{T}}dW_{\epsilon,H,\gamma}(s),
\end{equation}
where the integration is performed over the random measure $dW_{\epsilon,H,\gamma}$ defined in \ref{eq:DefdWRandom}. In this statistically stationary range, it is a zero-average process for any $H$ and $\epsilon>0$. Note by $X_{H,\gamma}(t)$ the process which marginals coincide with the limiting behavior of those of $X_{\epsilon,H,\gamma}(t)$ as $\epsilon\to 0$, and call it a Multifractal fractional Ornstein-Uhlenbeck (MfOU) process. For any $H\in]0,1[$, independently of $\gamma$, the covariance of $X_{H,\gamma}$ coincides with the one of the underlying fractional Ornstein-Uhlenbeck (fOU) process $X_{H,0}$. In particular, its variance is finite (\ref{eq:VarfOU}) and its covariance can be derived in an exact fashion (\ref{eq:RemExprCorrFOU1}).

More generally, considering the integer $n\ge 1$, and for any $H\in]0,1[$ and $\gamma^2<\min(1/4,H/(n-1))$, the MfOU process has a finite moment of order $2n$, that is
\begin{equation}
\E \left(X_{H,\gamma}^{2n} \right)< \infty,
\end{equation}
and we can obtain their exact expression as multiple integrals (see Section \ref{sec:Proofs} devoted to proofs, and in particular Eqs. \ref{eq:MomX2nHge12Final}, \ref{eq:MomX2nH=12Final} and \ref{eq:MomX2nHle12Final}). Furthermore, the MfOU process exhibits at small scales a multifractal behavior, as it can be characterized in a simple fashion by the behavior at small scales of the moments of its increments. We have, for an integer $n\ge 1$, $H\in]0,1[$ and $\gamma^2<\min(1/4,H/(n-1))$,
\begin{align}
\E\left( \delta_\tau X_{H,\gamma}\right)^{2n}\build{\sim}_{\tau\to 0^+}^{} c_{H,\gamma,2n}\frac{(2n)!}{2^nn!}\,\left(\frac{\tau}{T}\right)^{2nH-2n(n-1)\gamma^2},\label{eq:S2nMFOU}
\end{align}
where the multiplicative constant $c_{H,\gamma,2n}$ reads
\begin{align}\label{eq:CHGamma2n}
c_{H,\gamma,2n}&=T^{2nH}e^{2n(n-1)\gamma^2g(0)}\times\\
&\int_{\{u_i\}_{1\le i\le n}\in\R^{n}} \prod_{i=1}^n\left[(1-u_i)^{H-\frac{1}{2}}1_{u_i\le 1}-(-u_i)^{H-\frac{1}{2}}1_{u_i\le 0}\right]^2 \prod_{i<j=1}^n\frac{1}{|u_i-u_j|^{4\gamma^2}}\prod_{i=1}^{n}du_i,\notag
\end{align}
with $g(0)$ given in \ref{eq:G(0)}.

\end{proposition}
\smallskip

\begin{remark}Because of the independence of the multiplicative chaos $M_{\epsilon,\gamma}$ and the Wiener process $W$ entering in the definition of the random measure $dW_{\epsilon,H,\gamma}$ (\ref{eq:DefdWRandom}), the MfOU process $X_{H,\gamma}$ has a symmetric probability distribution. In consequence, all odd moments of $X_{H,\gamma}$ and its increments vanish.
\end{remark}

\smallskip

\begin{remark}Let us examine the values of the factors  $c_{H,\gamma,2n}$ in simple cases. First, their expression (\ref{eq:CHGamma2n}) is understood when $n=1$ as being
\begin{align}\label{eq:CHGamma2}
c_{H,\gamma,2}=c_{H,0,2}&=T^{2H}\int_{u\in\R} \left[(1-u)^{H-\frac{1}{2}}1_{u\le 1}-(-u)^{H-\frac{1}{2}}1_{u\le 0}\right]^2 du\\
&=\frac{T^{2H}}{\sin(\pi H)}\frac{\left[ \Gamma\left( H+\frac{1}{2}\right)\right]^2}{\Gamma\left( 2H+1\right)},\notag
\end{align}
which is independent of $\gamma$, and coincides with the one of the fOU process (its statistical properties are recalled in Remark \ref{Rem:fOU}), as it is claimed in Proposition \ref{Prop:MfOU}. Also, it is straightforward to get these factors for the non-intermittent $\gamma=0$ case since the MfOU process corresponds to the fOU process, and is thus Gaussian. It is then clear, looking at the expression provided in \ref{eq:CHGamma2n}, that
\begin{align}
c_{H,0,2n}=\left(c_{H,0,2}\right)^n,
\end{align}
as it is expected for a Gaussian process. The last situation for which the expression of these factors $c_{H,\gamma,2n}$ (\ref{eq:CHGamma2n}) gets simple is when $H=1/2$. In this case, we obtain
\begin{align}\label{eq:C12Gamma2n}
c_{1/2,\gamma,2n}&=T^{n}e^{2n(n-1)\gamma^2g(0)}
\int_{\{u_i\}_{1\le i\le n}\in[0,1]^{n}} \prod_{i<j=1}^n\frac{1}{|u_i-u_j|^{4\gamma^2}}\prod_{i=1}^{n}du_i,
\end{align}
that can be explicitly computed for the fourth-order moment (i.e. $n=2$) as 
\begin{align}\label{eq:C12Gamma4}
c_{1/2,\gamma,4}&=T^{2}e^{4\gamma^2g(0)}
\int_{\{u_1,u_2\}\in[0,1]^{2}} \frac{1}{|u_1-u_2|^{4\gamma^2}}du_1du_2\\
&=T^{2}e^{4\gamma^2g(0)}\frac{1}{(1-4\gamma^2)(1-2\gamma^2)},\notag
\end{align}
as it was already derived in \cite{VigFri20}. Exact values of $c_{H,\gamma,2n}$ for any appropriate $H$ and $\gamma$ are not known, but can be evaluated while performing a numerical estimation of the integrals entering in \ref{eq:CHGamma2n}.
\end{remark}

\smallskip

\begin{remark} Notice that the characterization of the sample path properties of $X_{H,\gamma}$ and its multifractal nature, as it is announced in Proposition \ref{Prop:MfOU}, using only its even moments of integer orders, and those of its increments, is rather simple and incomplete. A more complete analysis of the statistical properties of the MfOU process, in the spirit of \cite{CheGar19}, remains to be done.
\end{remark}

\section{Numerical synthesis and statistical analysis}\label{Sec:Nums}

 \subsection{Numerical synthesis}
The purpose of this Section is to propose a simple and efficient numerical algorithm that generates a trajectory of the statistically stationary MfOU process $X_{\epsilon,H,\gamma}(t) $ (\ref{eq:PropMfOU}), using similar ideas as they are developed in \cite{VigFri20}. As it is considered in the theoretical framework developed in Proposition \ref{Prop:MfOU}, trajectories are given for a finite $\epsilon>0$, and for $H\in]0,1[$ and a certain value of the parameter $\gamma$, the values of $H$ and $\gamma$ being possibly conditioned on the existence of moments of a certain order as $\epsilon\to 0$. Alternatively, we could solve the underlying dynamics itself (\ref{eq:SDEintro}) and thus access the transient towards the statistically regime using a time-marching discretization, as it is proposed in \cite{Che17,PerMor18}. In any case, we need to come up with a numerical approximation of the underlying random measure  $dW_{\epsilon,H,\gamma}(t)$ (\ref{eq:DefdWRandom}) which can be seen as a convolution of some random measure, that we will focus on in the sequel, with a causal and deterministic kernel, making the overall dynamics (\ref{eq:SDEintro}) non-Markovian. It is then tempting to use the Discrete Fourier transform (DFT) to design a numerical approximation $\widehat{dW}_{\epsilon,H,\gamma}[t]$ of the continuous measure $dW_{\epsilon,H,\gamma}(t)$, and incidentally, work in a framework for $T_{\text{\tiny{tot}}}$-periodic functions.

The causal and deterministic kernel involved in the definition of $dW_{\epsilon,H,\gamma}(t)$ reads
\begin{equation}\label{eq:DefdKernelLinearCombiNumerics}
h_{\epsilon,H}(t)= \mathfrak{h}_{\epsilon,H}(t)+\epsilon^{H-\frac{1}{2}}\delta(t),
\end{equation}
where
\begin{equation}\label{eq:DefdKernelhfrakNumerics}
\mathfrak{h}_{\epsilon,H}(t)= \left(H-\frac{1}{2}\right)(t+\epsilon)^{H-\frac{3}{2}}1_{t\ge 0}.
\end{equation}
In this continuous setting, the function $h_{\epsilon,H}$ (\ref{eq:DefdKernelLinearCombiNumerics}) behaves differently depending of the value of the Hurst exponent $H$, as it is classically interpreted when considering fractional Gaussian noises \cite{ManVan68}. For $H<1/2$, the function $h_{\epsilon,H}$ is negative outside the origin, and is integrable over $t\in\R$, such that its integral vanishes. On the opposite, when $H>1/2$, its integral diverges, whereas it is equal to unity for $H=1/2$. Notice also that, in the vicinity of the origin, $h_{\epsilon,H}$ is distributional, with a weight that depends on $\epsilon$.

From a numerical point of view, over a finite set $t\in[-T_{\text{\tiny{tot}}}/2,T_{\text{\tiny{tot}}}/2]$, we need to make a compromise between enforcing integral properties of $h_{\epsilon,H}$ (\ref{eq:DefdKernelLinearCombiNumerics}), especially for $H<1/2$, and respecting the strength of its distributional nature in the vicinity of the origin. The most natural way to build a numerical estimation $\widehat{h}_{\epsilon,H}$ of the function $h_{\epsilon,H}$ (\ref{eq:DefdKernelLinearCombiNumerics}) is obtained in the Fourier space, such that
\begin{align}
\widehat{h}_{\epsilon,H}[t]=\text{DFT}^{-1}\left\{ \text{DFT}\{\mathfrak{h}_{\epsilon,H}\}[\omega]+\epsilon^{H-\frac{1}{2}}\right\}[t],\label{eq:EstimhHin01} 
\end{align}
where the function $\mathfrak{h}_{\epsilon,H}$ is defined in \ref{eq:DefdKernelhfrakNumerics}, properly periodized over $[-T_{\text{\tiny{tot}}}/2,T_{\text{\tiny{tot}}}/2]$, and the frequency $\omega$ belonging to the discrete set  $[-1/(2T_{\text{\tiny{tot}}}),1/(2T_{\text{\tiny{tot}}})]$. It is clear that integral properties of $h_{\epsilon,H}$ (\ref{eq:DefdKernelLinearCombiNumerics}) are not exactly fulfilled by the proposed estimation $\widehat{h}_{\epsilon,H}$. In particular, the integral over $[-T_{\text{\tiny{tot}}}/2,T_{\text{\tiny{tot}}}/2]$ is equal, up to numerical errors regarding the finiteness of the time step, to $(T_{\text{\tiny{tot}}}/2+\epsilon)^{H-1/2}$, which eventually goes towards the expected limits as $T_{\text{\tiny{tot}}}\to \infty$. In the same time, the estimation reproduces accurately the distributional nature of the function $h_{\epsilon,H}$ close to the origin.

An important numerical observation made in \cite{VigFri20} shows that constraining the estimation $\widehat{h}_{\epsilon,H}$ to fulfill the integral properties of $h_{\epsilon,H}$ is crucial for the development of logarithmic correlations for the particular case $H=0$ (\ref{eq:DivCovXtilde}), and brings numerical stability. Recall that for this boundary case $H=0$, very different behaviors are expected. In particular, whereas the variance of the fOU process (\ref{eq:fOUProcInt}), when $H\in]0,1[$, remains bounded as $\epsilon$ goes to zero (\ref{eq:VarfOU}), it is expected to diverge when $H=0$ (\ref{eq:IntegralXtildeDef}) in a logarithmic fashion (\ref{eq:DivVarXtilde}).
Because of this observation and for these reasons, we treat the $H=0$ case separately and propose the following estimator $\widehat{h}_{\epsilon,0}$ of the function $h_{\epsilon,0}$ as
\begin{align}
\widehat{h}_{\epsilon,0}[t]=\text{DFT}^{-1}\left\{ \text{DFT}\{\mathfrak{h}_{\epsilon,0}\}[\omega]-\text{DFT}\{\mathfrak{h}_{\epsilon,0}\}[0]\right\}[t] \label{eq:EstimhH0}.
\end{align}
It is now clear that the vanishing integral properties of $h_{\epsilon,0}$  is exactly fulfilled by its estimation $\widehat{h}_{\epsilon,0}$ (\ref{eq:EstimhH0}), whereas an unavoidable error is made on the strength of the Dirac distribution near the origin.

Given these precisions, call $N$ the number of collocation points and $\Delta t=T_{\text{\tiny{tot}}}/N$ the time step, we propose the following numerical estimation $\widehat{\widetilde{X}}_{\epsilon}[t]$ of the regularized fOU process of vanishing Hurst exponent $\widetilde{X}_{\epsilon}(t)$ (\ref{eq:IntegralXtildeDef}) as
\begin{align}\label{eq:EstimXtildeHat}
\widehat{\widetilde{X}}_{\epsilon}[t]=\Delta t\;\text{DFT}^{-1}\left\{ \text{DFT}\left\{e^{-t/T}1_{t\ge 0}\right\}\text{DFT}\left\{\widehat{h}_{\epsilon,0}\right\}[\omega]\text{DFT}\left\{\widehat{d\widetilde{W}}\right\}\right\}[t],
\end{align}
where $\widehat{h}_{\epsilon,0}$ is given in \ref{eq:EstimhH0}, and $\widehat{d\widetilde{W}}$ a collection of $N$ independent instances of a zero average Gaussian random variable of variance $\Delta t$. To prevent from additional aliasing errors, consider long trajectories such that $T_{\text{\tiny{tot}}} \gg T$, where $T$ is the characteristic large scale of the OU kernel entering in \ref{eq:EstimXtildeHat}. Accordingly, an estimation $\widehat{M}_{\epsilon,\gamma}[t]$ of the random positive weight $M_{\epsilon,\gamma}(t)$ (\ref{eq:MCIntro}) is obtained while exponentiating the estimation $\widehat{\widetilde{X}}_{\epsilon}$ (\ref{eq:EstimXtildeHat}), such that
\begin{equation}\label{eq:EstimMC}
\widehat{M}_{\epsilon,\gamma}[t] = e^{\gamma \widehat{\widetilde{X}}_{\epsilon}[t]-\gamma^2\widehat{\E}\left(\left. \widehat{\widetilde{X}}_{\epsilon}\right.^2\right)},
\end{equation}
where enter the empirical expectation $\widehat{\E}$ of the square of the discrete process  $\widehat{\widetilde{X}}_{\epsilon}$. We could have used alternatively the exact expression of this expectation, that can be easily derived even for a given finite $\epsilon>0$ (see A13 of \cite{VigFri20}), which eventually gives undistinguishable results using the parameters of subsequent simulations.

Finally, the estimator $\widehat{X}_{\epsilon,H,\gamma}[t]$ of the statistically stationary MfOU process $X_{\epsilon,H,\gamma}(t) $ (\ref{eq:PropMfOU}) reads
\begin{align}\label{eq:EstimfOUProc}
\widehat{X}_{\epsilon,H,\gamma}[t]=\Delta t\;\text{DFT}^{-1}\left\{ \text{DFT}\left\{e^{-t/T}1_{t\ge 0}\right\}\text{DFT}\left\{\widehat{h}_{\epsilon,H}\right\}[\omega]\text{DFT}\left\{\widehat{M}_{\epsilon,\gamma}\widehat{dW}\right\}\right\}[t],
\end{align}
where $\widehat{h}_{\epsilon,H}$ is given in \ref{eq:EstimhHin01}, $\widehat{M}_{\epsilon,\gamma}$ in \ref{eq:EstimMC} and $\widehat{dW}$ a collection of $N$ independent instances of a zero average Gaussian random variable of variance $\Delta t$, independent of those of $\widehat{d\widetilde{W}}$.

\subsection{Statistical analysis}

We generate numerical instances of the discrete random process $\widehat{X}_{\epsilon,H,\gamma}$ (\ref{eq:EstimfOUProc}) taking without loss of generality $T_{\text{\tiny{tot}}}=1$. We consider $N=2^{30}$, and three different values for $H$, namely $H=1/3, \; 1/2$ and $2/3$, and three different values for $\gamma$, such that $\gamma^2=0, \;0.02$ and  $0.04$. For each couple of values $(H, \gamma^2)$, we attribute the values $T=T_{\text{\tiny{tot}}}/2^{10}$ for the characteristic large scale of the OU kernel entering in \ref{eq:EstimfOUProc}, and $\epsilon=4\Delta t$ concerning the regularizing scale, and generate ten independent trajectories.

\begin{figure}[t]
\begin{center}
\includegraphics[width=.9\linewidth]{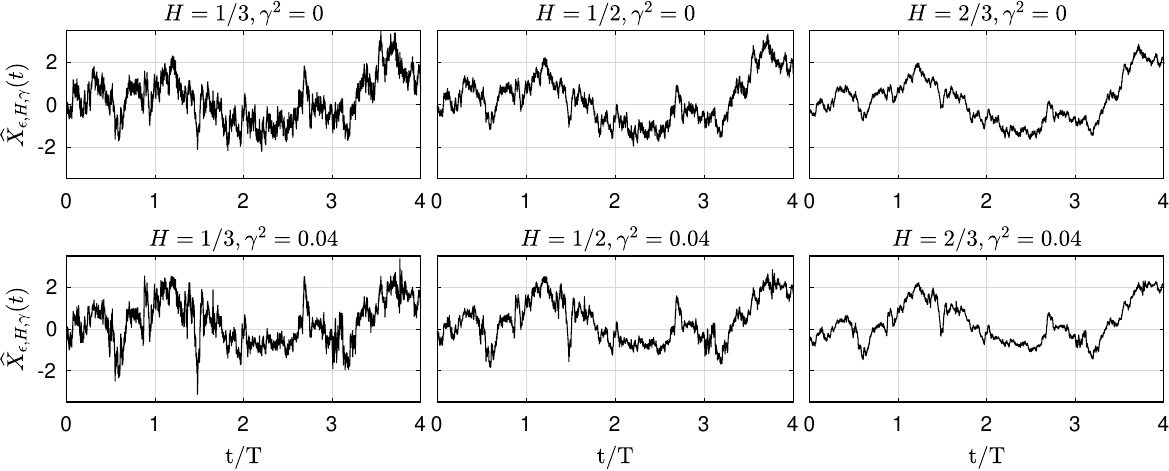}
\end{center}
\caption{An instance of a MfOU trajectory $\widehat{X}_{\epsilon,H,\gamma}$ (\ref{eq:EstimfOUProc}), for increasing values of the Hurst parameter $H$ (from left to right: $H =1/3$,  $1/2$ and $2/3$), and increasing values of the intermittent parameter (from top to bottom, $\gamma^2=0$ and $\gamma^2=0.04$). See the text for values of additional numerical parameters.}
\label{fig:traj}
\end{figure}

Figure \ref{fig:traj} shows an instance of the random process obtained $\widehat{X}_{\epsilon,H,\gamma}$ (\ref{eq:EstimfOUProc}) over several time scales $T$ and for different values of the couple $(H,\gamma)$, with the same instance of  the underlying independent white noises $\widehat{dW}$  and $\widehat{d\widetilde{W}}$. As $H$ increases, from left to right, the signal becomes less and less rough. Going from the Gaussian, i.e. non-intermittent, situation using $\gamma^2=0$ (top line), to its multifractal version using $\gamma^2=0.04$ (bottom line), we can see the appearance of spotty fluctuations, associated to events of larger time variation, which are reminiscent of the fluctuations of the random weight $\widehat{M}_{\epsilon,\gamma}$ (\ref{eq:EstimMC}).

\begin{figure}[t]
\begin{center}
\includegraphics[width=.9\linewidth]{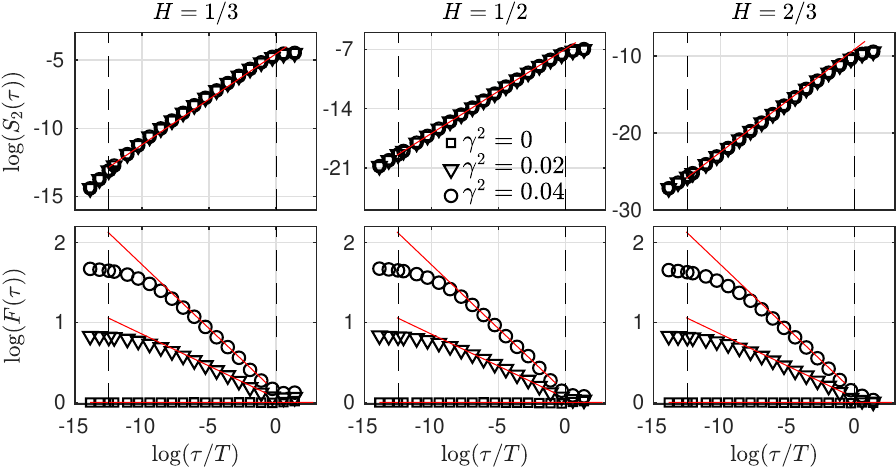}
\end{center}
\caption{Second-order moment $S_2(\tau)$ (\ref{eq:SnEmpirical}) (top line) and flatness $F(\tau)$ (\ref{eq:DefFlat}) (bottom line) of the increments of the MFOU process for different  values of $H$ (from left to right $H=1/3$, $1/2$ and $2/3$) and different values of $\gamma^2$ 
($\gamma^2=0$ ($\star$), $\gamma^2=0.02$ ($\triangledown$) and $\gamma^2=0.04$ ($\square$), in a logarithmic representation. The red lines correspond to the theoretical predictions concerning the asymptotic regime $\epsilon\to 0$ and given in \ref{eq:S2nMFOU}. We have superimposed vertical dashed lines to show the 
characteristic scales $\epsilon$ and $T$ used for the simulation.
}
\label{fig:moments}
\end{figure}

To investigate and quantify further the multiscale nature of the MfOU process $X_{\epsilon,H,\gamma}$,  we present a statistical analysis of some moments of its increments over a given scale $\tau$, and compare to the predictions made in Proposition \ref{Prop:MfOU} of Section \ref{Sec:StatResults}. More sophisticated analysis could be invoked to quantify the expected Multifractal phenomenon \cite{WenAbr07,ArnAud08,JafSeu19}, although we choose to keep the analysis simple and focus on the power-law behavior with $\tau$ of these moments. To do so, we compute the empirical moments of the increments $\delta_\tau \widehat{X}_{\epsilon,H,\gamma}[t]=\widehat{X}_{\epsilon,H,\gamma}[t+\tau]-\widehat{X}_{\epsilon,H,\gamma}[t]$ of our estimation $\widehat{X}_{\epsilon,H,\gamma}$ (\ref{eq:EstimfOUProc}), and average over ten independent trajectories. In particular, we display in Fig. \ref{fig:moments} the results of our statistical analysis using the so-called second order structure function $S_2(\tau)$, obtained taking $n=2$ in the expression
\begin{equation}\label{eq:SnEmpirical}
S_n(\tau)\equiv \widehat{\E}\left[\left( \delta_\tau \widehat{X}_{\epsilon,H,\gamma}\right)^n\right]\equiv \frac{\Delta t}{T_{\text{\tiny{tot}}}-\tau}\sum_{i=1}^{N-\tau/\Delta t} \left( \delta_\tau \widehat{X}_{\epsilon,H,\gamma}[t_i] \right)^n,
\end{equation}
 and the respective renormalized flatness factor 
 \begin{equation}\label{eq:DefFlat}
F(\tau) = \frac{S_4(\tau)}{3S_2^2(\tau)},
\end{equation}
which coincides with unity for Gaussian processes.

We display in Fig. \ref{fig:moments} (top) in a logarithmic fashion our statistical estimation of the second-order structure function $S_2(\tau)$ (\ref{eq:SnEmpirical}), after having checked that all represented scales have converged from a  statistical point of view. As it is expected, and stated in Proposition \ref{Prop:MfOU}, second-order statistical properties are not impacted by multifractal corrections, and are thus independent of the parameter $\gamma$: different symbols of Fig. \ref{fig:moments}, representing different values for $\gamma$, are gathered on a same curve. We superimpose, using a red solid line the power-law predicted in Proposition \ref{Prop:MfOU} and provided in Eqs. \ref{eq:S2nMFOU} and \ref{eq:CHGamma2n}, as far as second-order $n=2$ is concerned (see \ref{eq:CHGamma2} for an explicit expression at second-order). We can see that the prediction reproduces accurately the decay of the variance of increments as the scale $\tau$ decreases, for both scaling and amplitude (no additional fitting parameter is involved), in the inertial range $\epsilon\ll \tau \ll T$, as it is called in the literature of fluid turbulence \cite{Fri95}. For scales below $\epsilon$, we recover the scaling of the underlying OU process, that predicts $S_2(\tau)$ proportional to $\tau$. Above $T$, $S_2(\tau)$ saturates towards two times the variance of the process (\ref{eq:VarfOU}), as expected from a statistically stationary process. Actually, the whole range of available scales can be predicted and successfully compared to our estimations using the expression provided in \ref{eq:RemExprCorrFOU1}, properly generalized to a finite $\epsilon>0$, as it can be derived using the material developed in Refs. \cite{Che17,VigFri20}.

Let us now focus on the increments flatness $F(\tau)$ (\ref{eq:DefFlat}) that we display in Fig. \ref{fig:moments} (bottom). We can observe, as it is expected from a Gaussian process, $F(\tau)$ gets independent on scale $\tau$ when $\gamma=0$ (represented with $\star$), coinciding with $F(\tau)=1$, consistently with the predictions listed in Proposition \ref{Prop:MfOU}. When $\gamma$ differs from zero (using $\triangledown$ for $\gamma^2=0.02$ and $\square$ for $\gamma^2=0.04$), we observe a power-law behavior in the inertial range. Amplitude and exponent of this power-law are fully captured by our asymptotic predictions (\ref{eq:S2nMFOU}), obtained in the limit $\epsilon\to 0$ and for small scales, that read
\begin{equation}\label{eq:ExplicitFlatPred}
F(\tau) \build{\sim}_{\tau\to 0^+}^{}  \frac{c_{H,\gamma,4}}{c_{H,\gamma,2}^2}  \left( \frac{\tau}{T}\right)^{-4\gamma^2},
\end{equation}
where the coefficients $c_{H,\gamma,2n}$ are provided in \ref{eq:CHGamma2n}. Whereas we can get an explicit expression when $H=1/2$, as it is recalled in 
\ref{eq:C12Gamma4}, we can get their numerical value performing a double numerical integration. Doing so, we observe that the coefficient in front of the expected power-law (\ref{eq:ExplicitFlatPred}) depends weakly on $H$ (data not shown), and remains close to the exact value obtained when $H=1/2$ (
\ref{eq:C12Gamma4}) for the range of values that are explored in Fig. \ref{fig:moments}. We superimpose these predictions using red lines, and observe that they collapse with our statistical estimation of the flatness in the inertial range $\epsilon\ll \tau \ll T$. For scales below $\tau$, the flatness saturates towards a plateau that depends in particular on the precise value of $\epsilon$. For scales above $T$, flatness also saturates, towards a value that depends on chosen parameters. Overall, we can compare with great success our asymptotic theoretical predictions (\ref{eq:S2nMFOU}) with a statistical analysis of simulated trajectories.

\begin{figure}[t]
\begin{center}
\includegraphics[width=.9\linewidth]{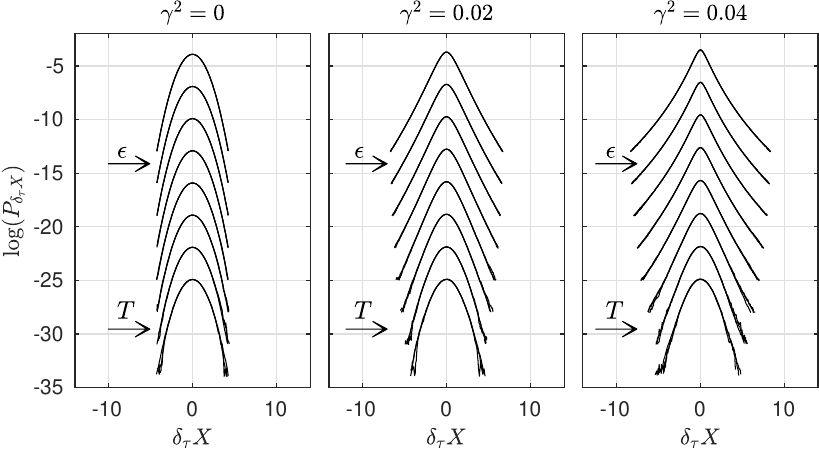}
\end{center}
\caption{Probability density function of increments, estimated while computing histograms, for different scales $\tau$. Shown scales $\tau$ are logarithmically spaced between $\epsilon/4$ and $4T$. For each scale, three curves are displayed corresponding to different values of the Hurst parameter $H$ ($H=1/3, 1/2$ and $2/3$), and are totally indistinguishable.}
\label{fig:PDFs}
\end{figure}

Finally, for the sake of completeness, we represent in Fig. \ref{fig:PDFs} the Probability Density Functions (PDFs) of increments, estimated as histograms, for different scales $\tau$ and for different values of $H$ and $\gamma$, as they are considered in Figs. \ref{fig:traj} and \ref{fig:moments}. For clarity, curves at different scales $\tau$ have been vertically shifted by an arbitraly factor, and are rescaled such that they are all of unit-variance. Firstly, we can observe that their shape depends indeed very weakly on the value of $H$, as we noticed while analyzing the values of the flatnesses (Fig. \ref{fig:moments}): we can barely distinguish them as $H$ varies. Also, we see that they remain Gaussian when $\gamma=0$ (left colomn), whereas they develop longer and longer tails as the scale decreases when $\gamma\ne 0$. Furthermore, all PDFs are symmetrical, as it is expected and recalled in Proposition \ref{Prop:MfOU}.

\section{Proofs}\label{sec:Proofs} 

Proving the statistical properties of the MfOU process $X_{H,\gamma}$, as they are stated in Proposition \ref{Prop:MfOU}, boils down to computing the moments of a random variable defined by a linear operation on the random measure $dW_{\epsilon,H,\gamma}$ (\ref{eq:DefdWRandom}), and then studying their limit when $\epsilon\to 0$. That is why we consider in the sequel the random variable
\begin{equation}\label{eq:MfOUConvolRV}
Y_{\epsilon,H,\gamma} = \int_{\R} f(s)dW_{\epsilon,H,\gamma}(s),
\end{equation}
where $f$ is an appropriate deterministic test function, say bounded and compactly supported for the sake of presentation. For the time being, assume furthermore that $f$ is differentiable over $\R$, note by $f'$ its derivative, an assumption that we will revisit and adapt when we will treat the particular case of the MfOU process $X_{H,\gamma}$ in Section \ref{Sec:PartCaseMFOU}. We recall here for convenience the expression of the random measure $dW_{\epsilon,H,\gamma}$ that enters in the definition of the random variable $Y_{\epsilon,H,\gamma}$ (\ref{eq:MfOUConvolRV}), as it was done in Definition \ref{Def:MFOU}. We have
\begin{equation}\label{eq:DefdWRandomProof}
dW_{\epsilon,H,\gamma}(s) = \omega_{\epsilon,H,\gamma}(s)ds + \epsilon^{H-\frac{1}{2}}M_{\epsilon,\gamma}(s)dW(s),
\end{equation}
where 
\begin{equation}\label{eq:OmegaHGammaProof}
\omega_{\epsilon,H,\gamma}(s) = \left(H-\frac{1}{2}\right)\int_{-\infty}^s(s-u+\epsilon)^{H-\frac{3}{2}}M_{\epsilon,\gamma}(u)dW(u),
\end{equation}
and $M_{\epsilon,\gamma}$ a regularized version of a Multiplicative Chaos (see Remark \ref{rem:DefMC}), normalized such that $\E M^2_{\epsilon,\gamma}=1$. Recall also that  $M_{\epsilon,\gamma}$ (\ref{eq:MCIntro}) is chosen independently of $dW$. It is then straightforward to see that $Y_{\epsilon,H,\gamma}$ is a zero-average random variable, and has more generally vanishing odd moments. In the sequel, we will find convenient to express the random measure $dW_{\epsilon,H,\gamma}$ (\ref{eq:DefdWRandomProof}) as a convolution product of the random measure $M_{\epsilon,\gamma}dW$, i.e.
\begin{equation}\label{eq:DefdWRandomProofLinearCombi}
dW_{\epsilon,H,\gamma}(s) =  ds\int_{\R} h_{\epsilon,H}(s-u)M_{\epsilon,\gamma}(u)dW(u),
\end{equation}
with the deterministic kernel
\begin{equation}\label{eq:DefdKernelLinearCombi}
h_{\epsilon,H}(s)= \left(H-\frac{1}{2}\right)(s+\epsilon)^{H-\frac{3}{2}}1_{s\ge 0}+\epsilon^{H-\frac{1}{2}}\delta(s),
\end{equation}
as it was already considered in Section \ref{Sec:Nums} devoted to numerics (see \ref{eq:DefdKernelLinearCombiNumerics}).

\subsection{Second-order statistical properties}\label{SecProof:SecondOS}

To prepare for subsequent calculations, first consider second-order statistical properties, and compute the variance of $Y_{\epsilon,H,\gamma}$, in a similar way, although more general, as it has been done in \cite{Che17}. We have
\begin{align}\label{eq:DefMom2Y}
\E\left[Y_{\epsilon,H,\gamma}^2\right] &= \int_{\R^2} f(s_1)f(s_2)\E\left[dW_{\epsilon,H,\gamma}(s_1)dW_{\epsilon,H,\gamma}(s_2)\right]\\
&= \int_{\R^3} f(s_1)f(s_2)h_{\epsilon,H}(s_1-u)h_{\epsilon,H}(s_2-u)\E\left[M^2_{\epsilon,\gamma}(u)\right]ds_1ds_2du\notag\\
&= \int_{\R} \mathcal L_{\epsilon,H}(u)du,\notag
\end{align}
where we have used that $M_{\epsilon,\gamma}$ is independent of $dW$ and normalized such that $\E (M_{\epsilon,\gamma}^2)=1$ (see Remark \ref{rem:DefMC}), and introduced the quantity
\begin{align}\label{eq:DefLEpsilon}
\mathcal L_{\epsilon,H}(u)=\int_{\R^2} f(s_1)f(s_2)h_{\epsilon,H}(s_1-u)h_{\epsilon,H}(s_1-u)ds_1ds_2.
\end{align}
This shows that there are no multifractal corrections on the second-order statistical properties of the random variable $Y_{\epsilon,H,\gamma}$, which coincide with those of $Y_{\epsilon,H,0}$. More explicitly, the deterministic function $\mathcal L_{\epsilon,H}$ (\ref{eq:DefLEpsilon}) reads
\begin{align}\label{eq:ExplicitLepsilonMom2Y}
&\mathcal L_{\epsilon,H}(u)= \left( H-\frac{1}{2}\right)^2  \int_{\R^2} f(s_1)f(s_2)(s_1-u+\epsilon)^{H-3/2}(s_2-u+\epsilon)^{H-3/2}1_{s_1\ge u}1_{s_2\ge u}ds_1ds_2\\
&+2\epsilon^{H-1/2} \left( H-\frac{1}{2}\right)\int_{\R}  f(s_1)f(u) (s_1-u+\epsilon)^{H-3/2}1_{s_1\ge u}ds_1 +\epsilon^{2H-1} f^2(u). \notag
\end{align}

For $H\in]1/2,1[$, we can see that, in the limit $\epsilon\to 0$, all integrals entering in $\mathcal L_{\epsilon,H}$ (\ref{eq:ExplicitLepsilonMom2Y}) involve singular kernels that are all locally integrable, thus only the first term of the RHS survives. Call this limit the function $\mathcal L_{H}(u)$. It is a continuous and bounded function, which is furthermore compactly supported if it is assumed so for the test function $f$. In this case, we can write
\begin{align}\label{eq:VarYHGE12final}
\E\left[Y_{H>1/2,\gamma}^2\right]&\equiv \lim_{\epsilon\to 0}\E\left[Y_{\epsilon,H,\gamma}^2\right]\equiv \int_{\R} \mathcal L_{H}(u)du\\
& = \left( H-\frac{1}{2}\right)^2  \int_{\R} \left( \int_{\R}f(s)(s-u)^{H-3/2}1_{s\ge u}ds\right)^2du.\notag
\end{align}
This last expression could be even further simplified, as it was done in \cite{Che17}, while introducing the correlation product, although, as we will see, this procedure is difficult to generalize for higher order moments. Notice that this expression makes perfect sense since all singularities are locally integrable, and $f$ is compactly supported.

For $H=1/2$, the situation is even simpler since the process $\omega_{\epsilon,H,\gamma}$ (\ref{eq:OmegaHGammaProof}) vanishes, and we obtain
\begin{align}
&\E\left[Y_{H=1/2,\gamma}^2\right]= \int_{\R} f^2(s)ds.
\end{align}

Take now $H\in]0,1/2[$. Going back to the expression of $\mathcal L_{\epsilon,H}$ (\ref{eq:ExplicitLepsilonMom2Y}), we can see that the first term in the RHS has no pointwise limit when $\epsilon\to 0$ since singularities are non locally integrable, neither do the two other terms. To extract finite and diverging contributions, split the integral of the first term of the RHS of \ref{eq:ExplicitLepsilonMom2Y} in two parts while considering the integration over the sets $s_1\le s_2$ and $s_1\ge s_2$, and notice that they give an equal contribution. Perform then an integration by parts over the dummy variable $s_1$ and obtain
\begin{align}\label{eq:SplitContOmegaInTwo}
  &2\left( H-\frac{1}{2}\right)^2  \int_{\R^2} f(s_1)f(s_2)(s_1-u+\epsilon)^{H-3/2}(s_2-u+\epsilon)^{H-3/2}1_{s_1\ge u}1_{s_2\ge s_1}ds_1ds_2\\
&=\mathcal I_{\epsilon,H}(u) + \mathcal J_{\epsilon,H}(u),\notag
\end{align}
where we have introduced the quantities
\begin{align}
\mathcal I_{\epsilon,H}(u)=\mathcal I_{1,H,\epsilon}(u)-\mathcal I_{2,H,\epsilon}(u)\notag
\end{align}
with
\begin{align}\label{eq:DefI1EpsilonHOM}
\mathcal I_{1,H,\epsilon}(u)&=2\left(H-\frac{1}{2}\right)^{}\int f^2(s)(s_{}-u+\epsilon)^{2H-2}1_{u\le s_{}}ds\\
&=-\epsilon^{2H-1}f^2(u)-2\int f(s)f'(s)(s_{}-u+\epsilon)^{2H-1}1_{u\le s_{}}ds,\notag \\
\mathcal I_{2,H,\epsilon}(u)&=2\left(H-\frac{1}{2}\right)^{}\epsilon^{H-\frac{1}{2}}\int f(u)f(s)(s_{}-u+\epsilon)^{H-\frac{3}{2}}1_{u\le s_{}}ds,\label{eq:DefI2EpsilonHOM}
\end{align}
and
\begin{align}\label{eq:DefJEpsilonHOM}
\mathcal J_{\epsilon,H}(u)=-2\left(H-\frac{1}{2}\right)^{}\int f'(s_1)f(s_2)(s_{1}-u+\epsilon)^{H-\frac{1}{2}}(s_{2}-u+\epsilon)^{H-\frac{3}{2}}1_{u\le s_{1}}1_{s_1\le s_{2}}ds_1ds_2.
\end{align}
Adding everything up, noticing that $\mathcal I_{2,H,\epsilon}(u)$ (\ref{eq:DefI2EpsilonHOM}) will cancel with the second term of the RHS of \ref{eq:ExplicitLepsilonMom2Y}, and that the contribution of order $\epsilon^{2H-1}$ of $\mathcal I_{1,H,\epsilon}(u)$ (\ref{eq:DefI1EpsilonHOM}) will cancel with the third one, we obtain
\begin{align}\label{eq:LUEPSHLE12}
\mathcal L_{\epsilon,H}(u)=\mathcal J_{\epsilon,H}(u)+\mathcal K_{\epsilon,H}(u),
\end{align}
where $\mathcal J_{\epsilon,H}$ is defined in \ref{eq:DefJEpsilonHOM} and
\begin{align}\label{eq:DefKEpsilonHOM}
\mathcal K_{\epsilon,H}(u)&=-2\int f(s)f'(s)(s_{}-u+\epsilon)^{2H-1}1_{u\le s_{}}ds\\
&=-2\int f(s+u)f'(s+u)(s+\epsilon)^{2H-1}1_{s\ge 0}ds.\notag
\end{align}
Whereas the pointwise limit as $\epsilon\to 0$ of the function $\mathcal L_{\epsilon,H}$ as written in \ref{eq:ExplicitLepsilonMom2Y} is not obvious when $H<1/2$, it becomes clear once written as in \ref{eq:LUEPSHLE12}. To see this, remark that the function $\mathcal K_{\epsilon,H}$ (\ref{eq:DefKEpsilonHOM}) is continuous, bounded and its integral over $u\in\R$ vanishes. Note its pointwise limit as the function $\mathcal K_{H} (u)$. Similarly, while taking the limit $\epsilon\to 0$, the strongest singularity entering in the expression of $\mathcal J_{\epsilon,H}$ (\ref{eq:DefJEpsilonHOM}) is obtained along the diagonal $s_1=s_2=s=u$ and is of order $(s-u)^{2H-2}$, which is locally integrable in $\R^2$ as long as $H>0$. Call then $\mathcal J_{H}(u)$ the pointwise limit of $\mathcal J_{\epsilon,H}(u)$.

We finally obtain, after taking the limit $\epsilon\to 0$,
\begin{align}\label{eq:VarYepsTo0HLE12}
&\E\left[Y_{H<1/2,\gamma}^2\right]=\int_{\R}\mathcal L_{H}(u)du = \int_{\R}\mathcal J_{H}(u)+\mathcal K_{H}(u)du = \int_{\R}\mathcal J_{H}(u)du\\
&=-2\left( H-\frac{1}{2}\right) \int_{\R^3} f'(s_1)f(s_2)(s_1-u)^{H-1/2}(s_2-u)^{H-3/2}1_{u\le s_{1}}1_{s_1\le s_{2}}duds_1ds_2.\notag
\end{align}

\subsection{Higher-order moments}\label{Sec:HOM}

Let us now develop and generalize our calculations, as they were done to derive the variance of the random variable $Y_{\epsilon,H,\gamma}$ (\ref{eq:MfOUConvolRV}), for the  moment of order $2n$, with $n$ and integer greater than unity. As mentioned in Proposition \ref{Prop:MfOU}, it is clear that its probability law is symmetrical, as a direct consequence of the independence of the white noise field $dW(t)_{t\in\R}$ and the random positive weight $M_{\epsilon,\gamma}(t)_{t\in\R}$ (\ref{eq:MCIntro}). Thus all moments of odd order vanish. 

Start with introducing the following quantity that encodes the correlation structure of the multiplicative chaos $M_{\epsilon,\gamma}$ (\ref{eq:MCIntro}) and write
\begin{align}\label{eq:CorrelMNpoints}
\mathcal M_{\epsilon,\gamma}(\{u_i\}_{1\le i\le n})=\E\left[ \prod_{i=1}^n M^2_{\epsilon,\gamma}(u_i) \right],
\end{align}
which converges as $\epsilon\to 0$ towards a locally integrable function in $\R^n$ when $\gamma^2<\min[1/4,1/(2(n-1))]$ (see Remark \ref{rem:DefMC}), and that reads
\begin{align}\label{eq:CorrelMNpointsEpsTo0}
\mathcal M_{\gamma}(\{u_i\})\equiv \lim_{\epsilon\to 0}\mathcal M_{\epsilon,\gamma}(\{u_i\})=T^{2n(n-1)\gamma^2} \prod_{i<j=1}^n\frac{e^{4\gamma^2g(u_i-u_j)}}{|u_i-u_j|_{+}^{4\gamma^2}}.
\end{align}

We have
\begin{align}\label{eq:Mom2nYEpsAllH}
&\E \left(Y_{\epsilon,H\in]0,1[,\gamma}^{2n} \right)=\int_{\{s_i\}\in\R^{2n}} \prod_{i=1}^{2n} f(s_i) \E\left[\prod_{i=1}^{2n} 
dW_{\epsilon,H,\gamma}(s_i)\right]\\
&=\int_{\{s_i,u_i\}\in\R^{4n}} \prod_{i=1}^{2n} f(s_i)h_{\epsilon,H}(s_i-u_i) \E\left[\prod_{i=1}^{2n} 
dW(u_i)\right]\E\left[\prod_{i=1}^{2n} M_{\epsilon,\gamma}(u_i)\right]  \prod_{i=1}^{2n} ds_idu_i    \notag \\
&=\frac{(2n)!}{2^nn!}\int_{\{s_i,u_i\}\in\R^{3n}}  \mathcal M_{\epsilon,\gamma}(\{u_i\})  \prod_{i=1}^{n} f(s_{2i-1}) f(s_{2i}) h_{\epsilon,H}(s_{2i-1}-u_i) h_{\epsilon,H}(s_{2i}-u_i)ds_{2i-1}ds_{2i}du_i    \notag \\
&=\frac{(2n)!}{2^nn!}\int_{\{u_i\}\in\R^{n}} \mathcal M_{\epsilon,\gamma}(\{u_i\}) \prod_{i=1}^{n} \mathcal L_{\epsilon,H}(u_i)  du_i,    \notag 
\end{align}
where the function $\mathcal L_{\epsilon,H}(u)$ is defined in \ref{eq:DefLEpsilon}. Recall from Section 
\ref{SecProof:SecondOS} that the function $\mathcal L_{\epsilon,H}(u)$ converges pointwise as $\epsilon\to 0$ towards a function $\mathcal L_{H}(u)$ which expression depends on the value of the parameter $H$. In particular, we obtained
\begin{align}\label{eq:ExpressLHge12}
\mathcal L_{H}(u) &\build{=}_{}^{H>1/2} \left( H-\frac{1}{2}\right)^2 \left( \int_{\R}f(s)(s-u)^{H-3/2}1_{s\ge u}ds\right)^2\\  \label{eq:ExpressLHequal12}
& \build{=}_{}^{H=1/2} f^2(u)\\
&\build{=}_{}^{H<1/2} \mathcal J_{H}(u) +\mathcal K_{H}(u), \label{eq:ExpressLHle12}
\end{align}
with
\begin{align}\label{eq:ExplicitJHU}
\mathcal J_{H}(u) &=-2\left(H-\frac{1}{2}\right)^{}\int_{\R^2} f'(s_1)f(s_2)(s_{1}-u)^{H-\frac{1}{2}}(s_{2}-u)^{H-\frac{3}{2}}1_{u\le s_{1}}1_{s_1\le s_{2}}ds_1ds_2,\\ \label{eq:ExplicitKHU}
\mathcal K_{H}(u) &=-2\int_{\R} f(s+u)f'(s+u)s^{2H-1}1_{s\ge 0}ds.
\end{align}
Remark that for any $H\in]0,1[$, the functions $\mathcal L_{H}$ and $f$ share similar supports, which ensures that a product on them will be integrable at large arguments. So justifying the existence of the limit
\begin{align}\label{eq:MomY2nFinal}
\E \left(Y_{H,\gamma}^{2n} \right)&\equiv \lim_{\epsilon\to 0}\E \left(Y_{\epsilon,H,\gamma}^{2n} \right)\\
&=\frac{(2n)!}{2^nn!}\int_{\{u_i\}\in\R^{n}} \mathcal M_{\gamma}(\{u_i\}) \prod_{i=1}^{n} \mathcal L_{H}(u_i)  du_i,\notag
\end{align}
where the pointwise limit $\mathcal M_{\gamma}$ of $\mathcal M_{\epsilon,\gamma}$ is provided in \ref{eq:CorrelMNpointsEpsTo0}, and in particular determining the range of parameters $H$ and $\gamma$ for which this limit makes sense, requires more assumptions on the test function $f$ than only continuity and boundedness. Instead of treating this question from a general point of view, we propose in the following Section to focus on the test function suggested by the MfOU process $X_{H,\gamma}$, and to justify its statistical properties as they are stated in Proposition \ref{Prop:MfOU}.

\subsection{Statistical properties of the MfOU process}\label{Sec:PartCaseMFOU}

\subsubsection{Variance and higher order moments of the process}

Let us now go back to our initial problem, that is obtaining the marginals of the process $X_{\epsilon,H,\gamma}$ (\ref{eq:PropMfOU}), and their limit as the regularizing scale $\epsilon$ goes to zero. We make use of the statistical stationarity of this process, and consider expectations at the time $t=0$. Concerning the moments of the process, consider the particular test function
\begin{equation}\label{eq:TestFuncVar}
f(s)=e^{\frac{s}{T}}1_{s\le 0},
\end{equation}
such that indeed $X_{\epsilon,H,\gamma}(t)$ (\ref{eq:PropMfOU}) at the given time $t=0$ can be written as
\begin{equation}\label{eq:MfOUtequal0}
X_{\epsilon,H,\gamma}(0) = \int_{\R} f(s)dW_{\epsilon,H,\gamma}(s).
\end{equation}
As we have seen in Section \ref{Sec:HOM}, for vanishing regularizing scale $\epsilon\to 0$, its moment of order $2n$ (\ref{eq:MomY2nFinal}) brings into play the function $\mathcal L_{H} (u)$ which explicit expression depends of the value of $H$ (Eqs. \ref{eq:ExpressLHge12}, \ref{eq:ExpressLHequal12} and \ref{eq:ExpressLHle12}). We have
\begin{align}\label{eq:ExpressLHge12MfOU}
\mathcal L_{H}(u) &\build{=}_{}^{H>1/2} 1_{u\le 0}\left( H-\frac{1}{2}\right)^2 \left( \int_{s=u}^0 e^{\frac{s}{T}}(s-u)^{H-3/2}ds\right)^2\\
&\build{\sim}_{u\to 0^-}^{}(-u)^{2H-1},\notag
\end{align}
which shows that $\mathcal L_{H}$ is a bounded function that goes to zero at the origin as fast as a power-law. This behavior in the vicinity of the origin weakens the singular behavior of the multiplicative chaos, and extends its range of integrability (\ref{eq:CorrelMNpointsEpsTo0}). Hence, the $2n$-order moment of the MfOU process is given by, for $H\in]1/2,1[$ and $\gamma^2<\min\left(1/4,H/(n-1)\right)$,
\begin{align}\label{eq:MomX2nHge12Final}
\E \left(X_{H>1/2,\gamma}^{2n} \right)=\frac{(2n)!}{2^nn!}\int_{\{u_i\}\in\R^{n}} \mathcal M_{\gamma}(\{u_i\}) \prod_{i=1}^{n} \mathcal L_{H}(u_i)  du_i,
\end{align}
where the function $\mathcal L_{H}$ is given in \ref{eq:ExpressLHge12MfOU} and $\mathcal M_{\gamma}$ in \ref{eq:CorrelMNpointsEpsTo0}. The case $H=1/2$ is rather simple since $\mathcal L_{H}(u)=f^2(u)=\exp(2u/T)1_{u\le 0}$, and we obtain for $\gamma^2<\min\left(1/4,1/(2(n-1))\right)$
\begin{align}\label{eq:MomX2nH=12Final}
\E \left(X_{H=1/2,\gamma}^{2n} \right)=\frac{(2n)!}{2^nn!}\int_{\{u_i\}\in\R^{n}} \mathcal M_{\gamma}(\{u_i\}) \prod_{i=1}^{n} f^2(u_i)  du_i.
\end{align}
Take now $H\in]0,1/2[$, and remark that, using the expression of the test function $f$ (\ref{eq:TestFuncVar}),
\begin{align}\label{eq:DerivTestFuncVar}
f'(s)=\frac{1}{T}f(s)-\delta(s),
\end{align}
which shows that the function $f$ is not differentiable on $\R$ in a classical sense, as expected from a causal kernel. For this reason, the expression of $\mathcal J_H$ (\ref{eq:ExplicitJHU}) where enters the derivative $f'$ is questionable. Nonetheless, we can argue, going back to the integration by parts procedure that takes place in the decomposition proposed in \ref{eq:SplitContOmegaInTwo}, that it makes sense to use the distributional derivative (\ref{eq:DerivTestFuncVar}) in the expression of the function $\mathcal J_H$ (\ref{eq:ExplicitJHU}) that enters the expression of $\mathcal L_H$ (\ref{eq:ExpressLHle12}). In this case, the Dirac distribution entering in the expression of $f'$ (\ref{eq:DerivTestFuncVar}) gives no contribution, and we obtain
\begin{align}\label{eq:ExplicitJHUOUKernel}
\mathcal J_{H}(u) &=-\frac{2}{T}\left(H-\frac{1}{2}\right)^{}\int_{\R^2} f(s_1)f(s_2)(s_{1}-u)^{H-\frac{1}{2}}(s_{2}-u)^{H-\frac{3}{2}}1_{u\le s_{1}}1_{s_1\le s_{2}}ds_1ds_2\\ 
&=-\frac{2}{T}\left(H-\frac{1}{2}\right)^{}1_{u\le 0}\int_{[u,0]^2} e^{\frac{s_1+s_2}{T}}(s_{1}-u)^{H-\frac{1}{2}}(s_{2}-u)^{H-\frac{3}{2}}1_{s_1\le s_{2}}ds_1ds_2.\notag
\end{align}
It is a positive, bounded and integrable function. To see this, we have for large negative arguments
\begin{align}
\mathcal J_{H}(u) \build{\sim}_{u\to -\infty}^{}-\frac{2}{T}\left(H-\frac{1}{2}\right)^{}(-u)^{2H-2}\int_{(\R^-)^2} e^{\frac{s_1+s_2}{T}}1_{s_1\le s_{2}}ds_1ds_2=-T\left(H-\frac{1}{2}\right)^{}(-u)^{2H-2}\notag,
\end{align}
which is indeed integrable. To see the behavior of $\mathcal J_{H}(u) $ at the origin, rescale the dummy variables $s_1$ and $s_2$ of \ref{eq:ExplicitJHUOUKernel} by $-u$ and get
\begin{align}
\mathcal J_{H}(u) & \build{\sim}_{u\to 0^-}^{}-\frac{2}{T}\left(H-\frac{1}{2}\right)^{}(-u)^{2H}\int_{[-1,0]^2} (s_{1}+1)^{H-\frac{1}{2}}(s_{2}+1)^{H-\frac{3}{2}}1_{s_1\le s_{2}}ds_1ds_2\notag\\
& \build{\sim}_{u\to 0^-}^{}-\frac{1}{T}\frac{H-\frac{1}{2}}{H\left(H+\frac{1}{2}\right)}(-u)^{2H}
.\notag
\end{align}
Whereas the function $\mathcal J_{H}(u)$ remains bounded, the function $\mathcal K_{H}(u)$ entering in the expression of $\mathcal L_{H}(u)$ (\ref{eq:ExpressLHle12}) will be singular at the origin, although integrable. To see this, notice first that we cannot make sense in a simple fashion of the expression of the function $\mathcal K_H$ (\ref{eq:ExplicitKHU}) while injecting the distributional derivative of $f$ (\ref{eq:DerivTestFuncVar}) in it. Instead, go back to the integration by parts procedure that takes place in the decomposition proposed in \ref{eq:SplitContOmegaInTwo}, and obtain, using the causal kernel $f$ (\ref{eq:TestFuncVar}),
\begin{align}\label{eq:ExplicitKHUAsympt}
\mathcal K_{H}(u) &=\lim_{\epsilon\to 0}1_{u\le 0}\left[ 2\left(H-\frac{1}{2}\right)\int_u^0 e^{2\frac{s}{T}}(s-u+\epsilon)^{2H-2}ds+\epsilon^{2H-1}e^{2\frac{u}{T}}\right]\\
&=\lim_{\epsilon\to 0}1_{u\le 0}\left[ (-u+\epsilon)^{2H-1}-\frac{2}{T}\int_u^0 e^{2\frac{s}{T}}(s-u+\epsilon)^{2H-1}ds\right]\notag\\
&=1_{u\le 0}\left[ (-u)^{2H-1}-\frac{2}{T}\int_u^0 e^{2\frac{s}{T}}(s-u)^{2H-1}ds\right]\notag.
\end{align}
This shows that $\mathcal K_{H}(u)$ is singular in the vicinity of the origin, such that
\begin{align}\label{eq:SingKHNear0}
\mathcal K_{H}(u) \build{\sim}_{u\to 0^-}^{}(-u)^{2H-1},
\end{align}
since the second term of the third equality of \ref{eq:ExplicitKHUAsympt} is bounded in $u$, and actually behaves proportionally to $(-u)^{2H}$ near the origin. To check the integrability at large arguments, write $\mathcal K_{H}(u)$ (\ref{eq:ExplicitKHUAsympt}) as
\begin{align}\label{eq:ExplicitKHUAsymptTowardInf}
\mathcal K_{H}(u) &=1_{u\le 0}\left[ \frac{2}{T}\int_u^0 e^{2\frac{s}{T}}\left[(-u)^{2H-1}-(s-u)^{2H-1}\right]ds +e^{2\frac{u}{T}}(-u)^{2H-1}\right]\\
&\build{\sim}_{u\to -\infty}^{}-\frac{2}{T}(2H-1)(-u)^{2H-2}\int_{-\infty}^0 s\;e^{2\frac{s}{T}}ds \notag\\
&\build{\sim}_{u\to -\infty}^{}\frac{T}{2}(2H-1)(-u)^{2H-2}, \notag
\end{align}
which is indeed integrable at large arguments. Actually, it can be shown that its integral over $u\in\R$ vanishes. To see this, before taking the limit $\epsilon\to 0$ in \ref{eq:ExplicitKHUAsympt}, observe that the integral vanishes for any $\epsilon>0$, which is also true in the limit.

Hence, the existence of the $2n$-order moment of the MfOU process is governed by the singular behaviors of $\mathcal K_{H}(u)$ (\ref{eq:SingKHNear0}) and $\mathcal M_{\gamma}$ (\ref{eq:CorrelMNpointsEpsTo0}), and is given by, for $H\in]0,1/2[$ and $\gamma^2<\min\left(1/4,H/(n-1)\right)$, 
\begin{align}\label{eq:MomX2nHle12Final}
\E \left(X_{H<1/2,\gamma}^{2n} \right)=\frac{(2n)!}{2^nn!}\int_{\{u_i\}\in\R^{n}} \mathcal M_{\gamma}(\{u_i\}) \prod_{i=1}^{n} \mathcal L_{H}(u_i)  du_i,
\end{align}
where the function $\mathcal L_{H}$ is given by $\mathcal L_{H}=\mathcal J_{H}+\mathcal K_{H}$, with $\mathcal J_{H}$ given in \ref{eq:ExplicitJHUOUKernel} and $\mathcal K_{H}$ in \ref{eq:ExplicitKHUAsympt}.

\subsubsection{Variance and higher order moments of the increments of the process}

Proceed now with the calculation of the moments of velocity increments, and determine in particular their behavior at small scales. To do so, consider the particular test function
\begin{align}\label{eq:KernelOUIncr}
f(s)\equiv f_\tau(s) &=e^{-\frac{\tau-s}{T}}1_{\tau-s\ge 0}-e^{-\frac{-s}{T}}1_{-s\ge 0}\\
&=e^{\frac{s}{T}}\left( e^{-\frac{\tau}{T}}1_{\tau-s\ge 0}-1_{-s\ge 0}\right).\notag
\end{align}
such that the increment $\delta_{\tau} X_{\epsilon,H,\gamma}(t)\equiv X_{\epsilon,H,\gamma}(t+\tau)-X_{\epsilon,H,\gamma}(t)$ of the MfOU process $X_{\epsilon,H,\gamma}$ (\ref{eq:PropMfOU}) at the given time $t=0$ can be written as
\begin{equation}\label{eq:IncrMfOUtequal0}
\delta_\tau X_{\epsilon,H,\gamma}(0) = \int_{\R} f(s)dW_{\epsilon,H,\gamma}(s).
\end{equation}
First consider $H>1/2$. In this case, make use of \ref{eq:ExpressLHge12} to get
\begin{align}\label{eq:ExpressLHge12IncrMfOU}
\mathcal L_{H}(u) &\build{=}_{}^{H>1/2} \left( H-\frac{1}{2}\right)^2 \left( \int_{\R} e^{\frac{s}{T}}\left( e^{-\frac{\tau}{T}}1_{\tau-s\ge 0}-1_{-s\ge 0}\right)(s-u)^{H-3/2}1_{s\ge u}ds\right)^2.
\end{align}
Rescaling by the positive quantity $\tau$, we get
\begin{align}
\mathcal L_{H}(\tau u) &\build{=}_{}^{H>1/2} \left( H-\frac{1}{2}\right)^2 \tau^{2H-1} \left( \int_{\R} e^{\frac{\tau s}{T}}\left( e^{-\frac{\tau}{T}}1_{s\le 1}-1_{s\le 0}\right)(s-u)^{H-3/2}1_{s\ge u}ds\right)^2\\
&\build{\sim}_{\tau\to 0^+}^{}\left( H-\frac{1}{2}\right)^2\tau^{2H-1} \left( \int_{\R} \left( 1_{s\le 1}-1_{s\le 0}\right)(s-u)^{H-3/2}1_{s\ge u}ds\right)^2\notag\\
&\build{\sim}_{\tau\to 0^+}^{}\tau^{2H-1} \left[(1-u)^{H-1/2}1_{u\le 1}-(-u)^{H-1/2}1_{u\le 0}\right]^2.\notag
\end{align}
Similarly for the term associated to the multiplicative chaos $\mathcal M_{\gamma}(\{u_i\})$ (\ref{eq:CorrelMNpointsEpsTo0}), we have, once its variables rescaled by the scale factor $\tau$,
\begin{align}\label{eq:CorrelMNpointsEpsTo0Rescaled}
\mathcal M_{\gamma}(\{\tau u_i\})&=T^{2n(n-1)\gamma^2} \prod_{i<j=1}^n\frac{e^{4\gamma^2g[\tau(u_i-u_j)]}}{|\tau(u_i-u_j)|_{+}^{4\gamma^2}}\\
&\build{\sim}_{\tau\to 0^+}^{}\left(\frac{\tau}{T}\right)^{-2n(n-1)\gamma^2} e^{2n(n-1)\gamma^2 g(0)}\prod_{i<j=1}^n\frac{1}{|u_i-u_j|^{4\gamma^2}},\notag
\end{align}
where $g(0)$ is given in \ref{eq:G(0)}, such that, for $H\in]1/2,1[$ and $\gamma^2<\min\left(1/4,H/(n-1)\right)$, we have
\begin{align}\label{eq:MomIncrX2nHge12Final}
\E \left(\delta_\tau X_{H>1/2,\gamma} \right)^{2n}&=\frac{(2n)!}{2^nn!}\int_{\{u_i\}\in\R^{n}} \mathcal M_{\gamma}(\{u_i\}) \prod_{i=1}^{n} \mathcal L_{H}(u_i)  du_i\\
&= \frac{(2n)!}{2^nn!}\tau^n\int_{\{u_i\}\in\R^{n}} \mathcal M_{\gamma}(\{\tau u_i\}) \prod_{i=1}^{n} \mathcal L_{H}(\tau u_i)  du_i\notag\\
&\build{\sim}_{\tau\to 0^+}^{} c_{H,\gamma,2n}\frac{(2n)!}{2^nn!}\,\left(\frac{\tau}{T}\right)^{2nH-2n(n-1)\gamma^2},
\end{align}
with  $c_{H,\gamma,2n}$ provided in \ref{eq:CHGamma2n}. This proves the power law announced in \ref{eq:S2nMFOU}. For $H=1/2$, since $\mathcal L_{H}$ is rather simple in this case (\ref{eq:ExpressLHequal12}), we have
\begin{align}
\mathcal L_{H}(\tau u) &\build{=}_{}^{H=1/2} e^{2\frac{\tau u}{T}}\left( e^{-\frac{\tau}{T}}1_{u\le 1}-1_{u\le 0}\right)^2\\
&\build{\sim}_{\tau\to 0^+}^{}\left( 1_{u\le 1}-1_{u\le 0}\right)^2,\notag
\end{align}
which justifies the expression of the constant $c_{1/2,\gamma,2n}$ (\ref{eq:CHGamma2n}) and of  the power law announced in \ref{eq:S2nMFOU} for $H=1/2$, with again $\gamma^2<\min\left(1/4,1/(2(n-1))\right)$.

Take now $H\in]0,1/2[$ and write the derivative of the kernel $f$ (\ref{eq:KernelOUIncr}) as
\begin{align}\label{eq:DerivKernelOUIncr}
f'(s) &=\frac{1}{T}f(s)-\delta(\tau-s)+\delta(s).
\end{align}
which makes sense as a distribution. Once again, it can be shown that we can inject safely this distribution into the expression of $\mathcal J_{H}(u)$ (\ref{eq:ExplicitJHU}), and obtain, while noticing that the Dirac function centered on $\tau$ entering in \ref{eq:DerivKernelOUIncr} has no contribution,
\begin{align}
\mathcal J_{H}(u)&=-\frac{2}{T}\left(H-\frac{1}{2}\right)^{}\int f(s_1)f(s_2)(s_{1}-u)^{H-\frac{1}{2}}(s_{2}-u)^{H-\frac{3}{2}}1_{u\le s_{1}}1_{s_1\le s_{2}}ds_1ds_2\\
&-2\left(H-\frac{1}{2}\right)^{}(-u)^{H-\frac{1}{2}}1_{u\le 0}\int f(s_2)(s_{2}-u)^{H-\frac{3}{2}}1_{0\le s_{2}}ds_2.\notag
\end{align}
Rescaling by the positive quantity $\tau$, we get
\begin{align}\label{eq:ExplicitJHUAsymptIncrRescale}
\mathcal J_{H}(\tau u)&=-\frac{2}{T}\left(H-\frac{1}{2}\right)^{}\tau^{2H}\int f(\tau s_1)f(\tau s_2)(s_{1}-u)^{H-\frac{1}{2}}(s_{2}-u)^{H-\frac{3}{2}}1_{u\le s_{1}}1_{s_1\le s_{2}}ds_1ds_2\\
&-2\left(H-\frac{1}{2}\right)^{}\tau^{2H-1}(-u)^{H-\frac{1}{2}}1_{u\le 0}\int f(\tau s_2)(s_{2}-u)^{H-\frac{3}{2}}1_{0\le s_{2}}ds_2\notag\\
&\build{\sim}_{\tau\to 0^+}^{} -2\left(H-\frac{1}{2}\right)^{}\tau^{2H-1}(-u)^{H-\frac{1}{2}}1_{u\le 0}\int (1_{s_2\le 1}-1_{s_2\le 0})(s_{2}-u)^{H-\frac{3}{2}}1_{0\le s_{2}}ds_2\notag\\
&\build{\sim}_{\tau\to 0^+}^{} -2\tau^{2H-1}(-u)^{H-\frac{1}{2}}1_{u\le 0}\left[(1-u)^{H-\frac{1}{2}}-(-u)^{H-\frac{1}{2}}\right].\notag
\end{align}
As encountered when deriving the variance of the process, we cannot make sense of the expression of $\mathcal K_{H}(u)$ (\ref{eq:ExplicitKHU}) when the derivative of $f$ (\ref{eq:DerivKernelOUIncr}) is distributional. Instead, write $\mathcal K_H$ as
\begin{align}\label{eq:ExplicitKHUAsymptIncr}
&\mathcal K_{H}(u) \equiv\lim_{\epsilon\to 0}\left[ 2\left(H-\frac{1}{2}\right)\int_\R f^2(s)(s-u+\epsilon)^{2H-2}1_{s\ge u}ds+\epsilon^{2H-1}f^2(u)\right]\\
&= (\tau-u)^{2H-1}1_{u\le \tau}+(-u)^{2H-1}\left( 1-2e^{-\frac{\tau}{T}}\right)1_{u\le 0}-\frac{2}{T}\int_\R f^2(s)(s-u)^{2H-1}1_{s\ge u}ds.\notag
\end{align}
Once rescaled by $\tau>0$, we obtain
\begin{align}\label{eq:ExplicitKHUAsymptIncrRescale}
\mathcal K_{H}(\tau u)&= \tau^{2H-1}\left[(1-u)^{2H-1}1_{u\le 1}+(-u)^{2H-1}\left( 1-2e^{-\frac{\tau}{T}}\right)1_{u\le 0}\right]\\
&-\tau^{2H}\frac{2}{T}\int_\R f^2(\tau s)(s-u)^{2H-1}1_{s\ge u}ds\notag\\
&\build{\sim}_{\tau\to 0^+}^{}	\tau^{2H-1}\left[(1-u)^{2H-1}1_{u\le 1}-(-u)^{2H-1}1_{u\le 0}\right].\notag
\end{align}
Finally, summing up the equivalents of $\mathcal J_{H}(\tau u)$ (\ref{eq:ExplicitJHUAsymptIncrRescale}) and $\mathcal K_{H}(\tau u)$ (\ref{eq:ExplicitKHUAsymptIncrRescale}) as $\tau \to 0$, we get
\begin{align}
\mathcal L_{H}(\tau u) \build{\sim}_{\tau\to 0^+}^{H<1/2} \tau^{2H-1} \left[(1-u)^{H-1/2}1_{u\le 1}-(-u)^{H-1/2}1_{u\le 0}\right]^2,\notag
\end{align}
which justifies the expression of the constant $c_{H,\gamma,2n}$ (\ref{eq:CHGamma2n}) and of  the power law announced in \ref{eq:S2nMFOU} for $H<1/2$, with again $\gamma^2<\min\left(1/4,H/(n-1)\right)$.

\smallskip 
\noindent \textbf{Acknowledgements.} We thank Ivan Nourdin for bringing to our knowledge the early work on a regularization procedure presented in \cite{AloMaz00}. L.C. is partially supported by the Simons Foundation Award ID: 651475.

\small

\end{document}